\documentclass[fractalfract,article,accept,moreauthors,pdftex]{my-mdpi} 

\firstpage{1} 
\pubvolume{xx}
\issuenum{1}
\articlenumber{5}
\pubyear{2020}
\copyrightyear{2020}

\externaleditor{Submitted: 19 May 2020; Revised: 14 July 2020; Accepted: 30 July 2020; Published: 1 August 2020} 

\history{Published Open Access in: \emph{Fractal Fract.} 2020, 4(3), Art.~38, 10~pp.
{\changeurlcolor{black}\href{https://doi.org/10.3390/fractalfract4030038}{\detokenize{https://doi.org/10.3390/fractalfract4030038}}}}

\pdfoutput=1

\setitemize{parsep=6pt,itemsep=0pt,leftmargin=*,labelsep=5.5mm}
\setenumerate{parsep=6pt,itemsep=0pt,leftmargin=*,labelsep=5.5mm}
\setlist[description]{itemsep=0mm}   

\Title{A Stochastic Fractional Calculus with Applications\\ to Variational Principles}

\Author{Houssine Zine $^{\dagger,\ddagger}$\orcidA{} 
and Delfim F. M. Torres *$^{,\ddagger}$\orcidB{}}

\AuthorNames{Houssine Zine and Delfim F. M. Torres}

\address[1]{Center for Research and Development in Mathematics and Applications (CIDMA),
Department of Mathematics, University of Aveiro, 3810-193 Aveiro, Portugal;
zinehoussine@ua.pt}

\corres{Correspondence: delfim@ua.pt}

\firstnote{This research is part of first author's Ph.D. project, 
which is carried out at the University of Aveiro under 
the Doctoral Program in Applied Mathematics
of Universities of Minho, Aveiro, and Porto (MAP).} 

\secondnote{These authors contributed equally to this work.}

\abstract{We introduce a stochastic fractional calculus.
As an application, we present a stochastic fractional calculus 
of variations, which generalizes the fractional calculus 
of variations to stochastic processes. A stochastic fractional 
Euler--Lagrange equation is obtained, extending those available 
in the literature for the classical, fractional, 
and stochastic calculus of variations. To illustrate our main 
theoretical result, we discuss two examples: one derived 
from quantum mechanics, the second validated 
by an adequate numerical simulation.}

\keyword{fractional derivatives and integrals; stochastic processes; calculus of variations} 

\MSC{26A33, 49K05, 60H10}

\begin{document}
	

\section{Introduction}

A stochastic calculus of variations, which generalizes 
the ordinary calculus of variations to stochastic processes, 
was introduced in 1981 by Yasue, generalizing the Euler--Lagrange equation 
and giving interesting applications to quantum mechanics \cite{R3}.
Recently, stochastic variational differential equations 
have been analyzed for modeling infectious diseases 
\cite{MR3808514,MR3872140}, and stochastic processes have shown to
be increasingly important in optimization \cite{MR3934431}.

In 1996, fifteen years after Yasue's pioneer work \cite{R3}, 
the theory of the calculus of variations evolved in order 
to include fractional operators and better describe 
non-conservative systems in mechanics \cite{R2}. The subject is
currently under strong development \cite{MR3822307}. We refer the
interested reader to the introductory book \cite{R2}
and to \cite{R1,baleanu2016numerical,baleanu2020fractional}
for numerical aspects on solving fractional Euler--Lagrange 
equations. For applications of fractional-order models and variational
principles in epidemics, biology, and medicine, see 
\cite{Ali2019319,BALEANU2020109705,MR4081589,Yousef2019297,MR3999702}
and references therein.

Given the importance of both stochastic and fractional calculi
of variations, it seems natural to join the two subjects.
That is the main goal of our current work, i.e., to introduce 
a stochastic-fractional calculus of variations. For that,
we start our work by introducing new definitions:  
left and right stochastic fractional derivatives and integrals
of Riemann--Liouville and Caputo types for stochastic processes 
of second order, as a deterministic function resulting from the 
intuitive action of the expectation, on which we can compute 
its fractional derivative several times to obtain additional 
results that generalize analogous classical relations.
Our definitions differ from those already available 
in the literature by the fact that they are applied on 
second order stochastic processes, whereas known definitions,
for example, those in \cite{MR3274137,MR2155620,MR1815336,MR2038026},  
are defined only for mean square continuous second order stochastic process, 
which is a short family of operators. Moreover, available results
in the literature have not used the expectation, which we claim to be 
more natural, easier to handle and estimate, when applied to 
fractional derivatives by different methods of approximation, 
like those developed and cited in \cite{R1}. More than different, 
our definitions are well posed and lead to numerous results 
generalizing those in the literature, like integration by parts 
and Euler--Lagrange variational equations. 

The paper is organized as follows. In Section~\ref{sec:2},
we introduce the new stochastic fractional operators.
Their fundamental properties are then given in Section~\ref{sec:3}.
In particular, we prove stochastic fractional formulas 
of integration by parts (see Lemma~\ref{lemmaIbyP}). 
Then, in Section~\ref{sec:4}, we consider the basic
problem of the stochastic fractional calculus of variations
and obtain the stochastic Riemann--Liouville and Caputo
fractional Euler--Lagrange equations (Theorems~\ref{thm:SFE-Leq}
and \ref{thm:SFE-Leq:C}, respectively). Section~\ref{sec:5} 
gives two illustrative examples. We end with Section~\ref{sec:6} 
on conclusions and future perspectives.
	

\section{The Stochastic Fractional Operators}
\label{sec:2}

Let $(\Omega,F,P)$ be a probabilistic space, where $\Omega$ is a nonempty set, 
$F$ is a $\sigma$-algebra of subsets of $\Omega$, and $P$ is a probability 
measure defined on $\Omega$. A mapping $X$ from an open time interval $I$ 
into the Hilbert space $H=L_2(\Omega,P)$ is a stochastic process of second order 
in $\mathbb{R}$. We introduce the stochastic fractional operators 
by composing the classical fractional operators with the expectation $E$.

In what follows, the classical fractional operators are denoted using standard
notations \cite{MR1347689}: ${_a}{D}{_t^\alpha}$~and ${_t}{D}{_b^\alpha}$
denote the left and right Riemann--Liouville fractional derivatives
of order $\alpha$; ${_a}{I}{_t^\alpha}$ and ${_t}{I}{_b^\alpha}$
the left and right Riemann--Liouville fractional integrals 
of order $\alpha$; while the left and right Caputo 
fractional derivatives of order $\alpha$ are denoted by
${_a^C}{D}{_t^\alpha}$ and ${_t^C}{D}{_b^\alpha}$, respectively.
The new stochastic operators add to the standard notations an 's'
for ``stochastic''.

\begin{Definition}[Stochastic fractional operators]
\label{def:SFO}
Let $X$ be a stochastic process on $[a,b] \subset I$, 
$\alpha>0$, $n=[\alpha]+1$, such that $E(X(t))\in AC^n([a,b]\rightarrow \mathbb{R})$  
with $AC$ the class of absolutely continuous functions. Then,
\begin{enumerate}[leftmargin=13 mm,labelsep=5.5mm]
\item[(D1)] the left stochastic Riemann--Liouville 
fractional derivative of order $\alpha$ is given by
\begin{eqnarray*}
{_a^s}{D}{_t^ \alpha X(t)} 
&=& {_a}{D}{_t^\alpha[E(X_t)]} \\
&=& \dfrac{1}{\Gamma (n-\alpha )}\left( \dfrac{d}{dt}\right)^n 
\int _a^t (t-\tau)^{n-1-\alpha} E(X_\tau) d\tau, \quad t>a;
\end{eqnarray*} 
\item[(D2)] the right stochastic Riemann--Liouville  
fractional derivative of order $\alpha$ by
\begin{eqnarray*}
{_t^s}{D}{_b^ \alpha X(t)} 
&=& {_t}{D}{_b^\alpha[E(X_t)]} \\
&=& \dfrac{1}{\Gamma (n-\alpha )}\left( \dfrac{-d}{dt}\right)^n 
\int _t^b (\tau-t)^{n-1-\alpha} E(X_\tau) d\tau, \quad t<b;
\end{eqnarray*}
\item[(D3)] the left stochastic Riemann--Liouville 
fractional integral of order $\alpha$ by
\begin{eqnarray*}
{_a^s}{I}{_t^ \alpha X(t)} 
&=& {_a}{I}{_t^\alpha[E(X_t)]}\\
&=& \dfrac{1}{\Gamma (\alpha )} 
\int _a^t (t-\tau)^{\alpha-1} E(X_\tau) d\tau, \quad t>a;
\end{eqnarray*} 
\item[(D4)] the right stochastic  
Riemann--Liouville fractional integral of order $\alpha$ by
\begin{eqnarray*}
{_t^s}{I}{_b^ \alpha X(t)} 
&=& {_t}{I}{_b^\alpha[E(X_t)]} \\
&=& \dfrac{1}{\Gamma (\alpha )} 
\int_t^b (\tau-t)^{\alpha-1} E(X_\tau) d\tau, 
\quad t<b;
\end{eqnarray*}

\item[(D5)] the left stochastic Caputo 
fractional derivative of order $\alpha$ by
\begin{eqnarray*}
{_a^{sC}}{D}{_t^ \alpha X(t)} 
&=& {_a^C}{D}{_t^\alpha[E(X_t)]} \\
&=& \dfrac{1}{\Gamma (n-\alpha )} 
\int _a^t (t-\tau)^{n-1-\alpha} E(X(\tau))^{(n)} d\tau; 
\quad t>a.
\end{eqnarray*}
\item[(D6)] and the right stochastic Caputo fractional derivative 
of order $\alpha$ by
\begin{eqnarray*}
{_t^{sC}}{D}{_b^ \alpha X(t)} 
&=& {_t^C}{D}{_b^\alpha[E(X_t)]}\\
&=& \dfrac{(-1)^{n}}{\Gamma (n-\alpha )} 
\int_t^b (\tau-t)^{n-1-\alpha}E(X(\tau))^{(n)} d\tau, 
\quad t<b.
\end{eqnarray*}
\end{enumerate}
\end{Definition}

\begin{Remark}
The stochastic processes $X(t)$ used along the manuscript 
can be of any type satisfying the announced conditions 
of existence of the novel stochastic fractional operators. 
For example, we can consider Levy processes as a particular case, 
provided one considers some intervals where $E(X(t))$ is sufficiently 
smooth \cite{MR3825051}.
\end{Remark}

As we shall prove in the following sections, the new 
stochastic fractional operators just introduced
provide a rich calculus with interesting applications.


\section{Fundamental Properties}
\label{sec:3}

Several properties of the classical fractional operators,
like boundedness or linearity, also hold true for their
stochastic counterparts.

\begin{Proposition}
\label{prop:01}
If $t\rightarrow E(X_t) \in L_1([a,b])$, then ${_a^s}{I}{_t^\alpha(X_t)}$ is bounded.
\end{Proposition}

\begin{proof}
The property follows easily from definition (D3):
$$
\left| {_a^s}{I}{_t^\alpha(X_t)} \right| 
= \left| \dfrac{1}{\Gamma (\alpha )} 
\int _a^t (t-\tau)^{\alpha-1} E(X_\tau) d\tau \right| 
\leq k \left\| E (X_t) \right\|_1,
$$
which shows the intended conclusion.
\end{proof}

\begin{Proposition}
\label{prop:02}
The left and right stochastic Riemann--Liouville 
and Caputo fractional operators given 
in Definition~\ref{def:SFO} are linear operators.
\end{Proposition}

\begin{proof}
Let $c$ and $d$ be real numbers and assume that 
${_a^s}{D}{_t^\alpha X_t}$ and ${_a^s}{D}{_t^\alpha Y_t}$ exist.
It is easy to see that ${_a^s}{D}{_t^\alpha (c\cdot X_t+d\cdot Y_t)}$ 
also exists. From Definition~\ref{def:SFO} and
by linearity of the expectation and the linearity 
of the classical/deterministic fractional derivative operator, we have
\begin{equation*}
\begin{split}
{_a^s}{D}{_t^\alpha (c\cdot X_t+d\cdot Y_t)}
&={_a}{D}{_t^\alpha E(c\cdot X_t+d\cdot Y_t)}\\
&=c\cdot {_a}{D}{_t^\alpha E(X_t)}+d\cdot {_a}{D}{_t^\alpha E(Y_t)}\\
&=c\cdot {_a^s}{D}{_t^\alpha(X_t)}+d\cdot {_a^s}{D}{_t^\alpha (Y_t)}.
\end{split}
\end{equation*}

The linearity of the other stochastic fractional operators 
is obtained in a similar manner.
\end{proof}

Our next proposition involves both stochastic and deterministic operators.
Let $\mathcal{O} \in \left\{D, I, {^{C}{D}} \right\}$.
Recall that if ${_a^s}{\mathcal{O}}{_t^\beta}$ is a left stochastic fractional
operator of order $\beta$, then ${_a}{\mathcal{O}}{_t^\beta}$ is the corresponding 
left classical/deterministic fractional operator of order $\beta$; 
similarly for right operators.

Note that the proofs of Propositions~\ref{prop:03} and \ref{prop:04} 
and Lemma~\ref{lemmaIbyP} are not hard to prove in the sense that they 
are based on well-known results available for deterministic fractional 
derivatives (observe that $E(X(t))$ is deterministic). 

\begin{Proposition}
\label{prop:03}
Assume that ${_a^s}{I}{_t^\beta X_t}$, ${_t^s}{I}{_b^\beta X_t}$,
${_a^s}{I}{_t^\alpha X_t}$, 
${_a}{D}{_t^\alpha \left[{_a^s}{I}{_t^\alpha X_t} \right] }$, 
${_a}{I}{_t^\alpha \left[{_a^s}{I}{_t^\beta X_t} \right] }$ 
and ${_t}{I}{_b^\alpha \left[{_t^s}{I}{_b^\beta X_t} \right] }$ exist. 
The~following relations hold:
\begin{gather*}
{_a}{I}{_t^\alpha \left[{_a^s}{I}{_t^\beta X_t} \right]}
={_a^s}{I}{_t^{\alpha+\beta}} X_t, \\
{_t}{I}{_b^\alpha \left[{_t^s}{I}{_b^\beta X_t} \right] }
={_t^s}{I}{_b^{\alpha+\beta}} X_t,\\
{_a}{D}{_t^\alpha \left[{_a^s}{I}{_t^\alpha X_t} \right] }
= E(X_t).
\end{gather*}
\end{Proposition}

\begin{proof}
Using Definition~\ref{def:SFO} and well-known properties 
of the deterministic Riemann--Liouville fractional operators 
\cite{ALMTOR}, one has
\begin{equation*}
\begin{split}
{_a}{I}{_t^\alpha \left[{_a^s}{I}{_t^\beta X_t} \right]}
&= {_a}{I}{_t^\alpha \left[{_a}{I}{_t^\beta E( X_t)} \right]}\\
&= {_a}{I}{_t^{\alpha+\beta}}E(X_t)\\
&={_a^s}{I}{_t^{\alpha+\beta}} X_t.  
\end{split}
\end{equation*}

The second and third equalities are easily proved in a similar manner.
\end{proof}

\begin{Proposition}
\label{prop:04}
Let $\alpha > 0$.
If $E(X_t)\in L_\infty(a,b)$, then
$$
{_a^C}{D}{_t^\alpha \left[{_a^s}{I}{_t^\alpha X_t} \right]}= E( X_t)
$$
and
$$
{_t^C}{D}{_b^\alpha \left[{_t^s}{I}{_b^\alpha X_t} \right]}= E( X_t).
$$
\end{Proposition}

\begin{proof}
Using Definition~\ref{def:SFO} and well-known properties 
of the deterministic Caputo fractional operators~\cite{ALMTOR}, we have
\begin{equation*}
\begin{split}
{_a^C}{D}{_t^\alpha \left[{_a^s}{I}{_t^\alpha X_t} \right]}
&= {_a^C}{D}{_t^\alpha \left[{_a}{I}{_t^\alpha E( X_t)} \right]}\\ 
&=  E( X_t). 
\end{split}
\end{equation*}

The second formula is shown with the same argument.
\end{proof}

Formulas of integration by parts play a fundamental role 
in the calculus of variations and optimal control 
\cite{MR3221831,MR4015001}. Here we make use 
of Lemma~\ref{lemmaIbyP} to prove in Section~\ref{sec:4} 
a stochastic fractional Euler--Lagrange 
necessary optimality condition.

\begin{Lemma}[Stochastic fractional formulas of integration by parts]
\label{lemmaIbyP}
Let $\alpha>0$, $p,q \geq 1$, and $\frac{1}{p}+\frac{1}{q}\leq 1+\alpha$ 
($p\neq 1$ and $q\neq 1$ in the case where $\frac{1}{p}+\frac{1}{q}=1+\alpha)$.
\begin{itemize}[leftmargin=10 mm,labelsep=5.5mm]
\item[(i)] If $E(X_t) \in L_p(a,b)$ 
and $E(Y_t) \in L_q(a,b)$ for every $t\in [a,b]$, then
$$
E\left( \int_a^b (X_t) {_a^s}{I}{_t^ \alpha} Y_t dt \right)
= E\left( \int_a^b (Y_t) {_t^s}{I}{_b^ \alpha} X_t dt \right).
$$
\item[(ii)] If $E(Y_t) \in {_t}{I}{_b^{\alpha}}(L_p)$ and 
$E(X_t) \in {_a}{I}{_t^{\alpha}}(L_q)$ for every $t\in [a,b]$, then
$$
E\left( \int_a^b (X_t) ({_a^s}{D}{_t^ \alpha} Y_t) dt \right)
= E\left( \int_a^b (Y_t) ({_t^s}{D}{_b^ \alpha} X_t) dt \right).
$$
\item[(iii)] For the stochastic Caputo fractional derivatives, one has 
$$ 
E \left[ \int_a^b (X_t) ({_a^{sC}}{D}{_t^ \alpha} Y_t) dt \right]
= E \left[ \int_a^b (Y_t) ({_t^s}{D}{_b^{\alpha}} X_t )dt \right]
+E \left[ ({_t^s}{I}{_b^{1-\alpha}} X_t)\cdot Y_t \right]_a^b
$$
and
$$E \left[ \int_a^b (X_t) ({_t^{sC}}{D}{_b^ \alpha} Y_t) dt \right]
= E \left[ \int_a^b (Y_t) ({_a^s}{D}{_t^{\alpha}} X_t )dt \right]
-E \left[ ({_a^s}{I}{_t^{1-\alpha}} X_t)\cdot Y_t \right]_a^b
$$ 
for $\alpha\in(0,1)$. 
\end{itemize}
\end{Lemma}

\begin{proof}
\begin{itemize}[leftmargin=10 mm,labelsep=5.5mm]
\item[(i)] We have
\end{itemize}
\begin{equation*}
\begin{split}
E\left( \int_a^b (X_t) {_a^s}{I}{_t^ \alpha} Y_t dt \right) 
&=  \int_a^b E\left( (X_t) {_a^s}{I}{_t^ \alpha} Y_t \right) dt 
\quad (\text{by Fubini--Tonelli's theorem})\\
&= \int_a^b E\left( (X_t) {_a}{I}{_t^ \alpha} E(Y_t) \right) dt 
\quad (by  \ (D_3)) \\
&=  \int_a^b E((X_t)) {_a}{I}{_t^ \alpha} E(Y_t) dt 
\quad (\text{the expectation is  deterministic})\\
&=  \int_a^b {_t}{I}{_b^ \alpha}E(X_t) \cdot E(Y_t) dt  \quad
(\text{by fractional integration by parts})\\
&= E\left(\int_a^b {_t^s}{I}{_b^ \alpha}(X_t) (Y_t) dt\right) 
\quad (\text{by Fubini--Tonelli's theorem}).
\end{split}
\end{equation*}
\begin{itemize}[leftmargin=10 mm,labelsep=5.5mm]
\item[(ii)]  With similar arguments as in item (i), we have
\end{itemize}
\begin{equation*}
\begin{split}
E\left( \int_a^b (X_t) {_a^s}{D}{_t^ \alpha} Y_t dt \right) 
&=  \int_a^b E\left( (X_t) {_a^s}{D}{_t^ \alpha} Y_t \right) dt\\ 
&= \int_a^b E\left( (X_t) {_a}{D}{_t^ \alpha} E(Y_t) \right) dt 
\quad (by  \ (D_1)) \\
&=  \int_a^b E((X_t)) {_a}{D}{_t^ \alpha} E(Y_t) dt \\
&=  \int_a^b {_t}{D}{_b^ \alpha}E(X_t) \cdot E(Y_t) dt \\
&= E\left(\int_a^b {_t^s}{D}{_b^ \alpha}(X_t) (Y_t) dt\right).
\end{split}
\end{equation*}
\begin{itemize}[leftmargin=10 mm,labelsep=5.5mm]
\item[(iii)]  By using Caputo's fractional integration by parts formula we obtain that 
\end{itemize}
\begin{equation*}
\begin{split}
E \left[ \int_a^b (X_t) ({_a^{sC}}{D}{_t^ \alpha} Y_t) dt \right]
&=  \int_a^b E \left[(X_t)\right] ({_a^{C}}{D}{_t^ \alpha}E\left[ (Y_t)\right]) dt \\
&=  \int_a^b ({_t}{D}{_b^ \alpha}E\left[ (X_t)\right]\cdot E \left[(Y_t)\right]) dt 
+ \left[ ({_t}{I}{_b^{1-\alpha}} E(X_t)\cdot E(Y_t) \right]_a^b\\
&= \int_a^b ({_t^{s}}{D}{_b^ \alpha}(X_t)\cdot  E \left[(Y_t)\right]) dt 
+ \left[ ({_t}{I}{_b^{1-\alpha}} E(X_t)\cdot E(Y_t) \right]_a^b\\
&= E \left[ \int_a^b ({_t^{s}}{D}{_b^ \alpha}(X_t)\cdot  (Y_t))dt \right] 
+E \left[ ({_t}{I}{_b^{1-\alpha}} E(X_t)\cdot (Y_t) \right]_a^b. 
\end{split}
\end{equation*}

The first equality of (iii) is proved. By using a similar argument and applying 
the integration by parts formula associated with the right Caputo 
fractional derivative \cite{ALMTOR}, we easily get the second equality of (iii). 
\end{proof}


\section{Stochastic Fractional Euler--Lagrange Equations}
\label{sec:4}

Let us denote by $C^1(I\rightarrow H)$ the set of second order 
stochastic processes $X$ such that the left and right 
stochastic Riemann--Liouville fractional derivatives 
of $X$ exist, endowed with the norm   
$$
\Vert X \Vert = \sup_{t\in I}\left(\Vert X(t) \Vert_H 
+\mid {_a^s}{D}{_t^\alpha}X(t) \mid 
+\mid {_t^s}{D}{_b^\alpha}X(t)\mid\right),
$$
where $\Vert \cdot \Vert_H$ is the norm of $H$. 
Let $L \in C^{1}(I\times H\times\mathbb{R}
\times \mathbb{R}\rightarrow\mathbb{R})$ 
and consider the following minimization problem:
\begin{equation}
\label{eq:F}
J[X]=E\left(\int_a^b 
L\left( t,X(t),{_a^s}{D}{_t^ \alpha} X(t),
{_t^s}{D}{_b^ \alpha} X(t)\right)dt\right)
\longrightarrow \min
\end{equation}
subject to the boundary conditions
\begin{equation}
\label{eq:BC}
E(X(a))=X_a, \quad E(X(b))=X_b,
\end{equation}
where $X$ verifies the above conditions and $L$ is a smooth function.
Taking into account the method used in \cite{R1} for the fractional setting, 
and according to stochastic fractional integration by parts 
given by our Lemma~\ref{lemmaIbyP}, we obtain the following necessary 
optimality condition for the fundamental problem~\eqref{eq:F}--\eqref{eq:BC}
of the stochastic fractional calculus of variations.

\begin{Theorem}[The stochastic Riemann--Liouville fractional Euler--Lagrange equation]
\label{thm:SFE-Leq}
If $J \in C^{1}(H \times \mathbb{R} 
\times \mathbb{R} \rightarrow \mathbb{R})$ 
and $X\in C^1(I\rightarrow H)$ is an 
$F$-adapted stochastic process on $[a,b]$ 
with $E(X(t))\in AC([a,b])$ that is a minimizer 
of \eqref{eq:F} subject to the fixed end points \eqref{eq:BC}, 
then $X$ satisfies the following stochastic 
fractional Euler--Lagrange equation:
$$
\frac{\partial L}{\partial X}+{_t^s}{D}{_b^ \alpha}\left[ 
\dfrac{\partial L}{\partial{_a^s}{D}{_t^ \alpha}}\right]
+{_a^s}{D}{_t^ \alpha}\left[ 
\dfrac{\partial L}{\partial{_t^s}{D}{_b^ \alpha}}\right]=0.
$$
\end{Theorem}

\begin{proof}
We have 
$$
J[X]=E\left(\int_a^b L(t,X(t),{_a^s}{D}{_t^ \alpha} X(t),
{_t^s}{D}{_b^ \alpha} X(t)dt\right).
$$

Assume that $X^*$ is the optimal solution
of problem \eqref{eq:F}--\eqref{eq:BC}. Set
$$
X=X^*+\varepsilon \eta,
$$
where $\eta$ is an $F$-adapted stochastic process 
on $[a,b]$ in $C^1(I\rightarrow H)$. By linearity 
of the stochastic fractional derivatives 
(Proposition~\ref{prop:02}), we get
$$
{_a^s}{D}{_t^ \alpha} X={_a^s}{D}{_t^ \alpha} X^*
+\varepsilon\left( {_a^s}{D}{_t^ \alpha} \eta\right)
$$
and
$$
{_t^s}{D}{_b^ \alpha} X
={_t^s}{D}{_b^ \alpha} X^*
+\varepsilon\left( {_t^s}{D}{_b^ \alpha} \eta\right).
$$

Consider now the following function:
$$ 
j(\varepsilon)=E\left(\int_a^b L\left(t,X^*
+\varepsilon \eta,{_a^s}{D}{_t^ \alpha} X^*
+\varepsilon\left( {_a^s}{D}{_t^ \alpha} \eta\right),
{_t^s}{D}{_b^ \alpha} X^*
+\varepsilon\left( {_t^s}{D}{_b^ \alpha} \eta \right)\right) dt\right).
$$

We deduce, by the chain rule, that
$$
\dfrac{d}{dt}j(\varepsilon)\mid_{\varepsilon=0}
=E\left( \int_a^b( \partial_{2}L \cdot \eta
+\partial_{3}L \cdot {_a^s}{D}{_t^ \alpha}\eta
+\partial_{4}L \cdot {_t^s}{D}{_b^ \alpha}\eta)dt\right) =0,
$$
where $\partial_i L$ denotes the partial derivative of the Lagrangian
$L$ with respect to its $i$th argument. Using Lemma~\ref{lemmaIbyP} 
of stochastic fractional integration by parts, we obtain
$$
E\left( \int_a^b\left(  \partial_{2}L \cdot \eta+{_t^s}{D}{_b^ \alpha}(\partial_{3}L)
\cdot \eta+{_a^s}{D}{_t^ \alpha}(\partial_{4}L) \cdot \eta\right) dt\right) =0.
$$

We claim that if $Y$ is a stochastic process with continuous paths 
of second order such that
$$
E\left[\int_a^b Y(t)\cdot \eta(t)dt\right]=0
$$
for any stochastic process with continuous paths $\eta$, then
$$
Y=0 \quad \text{almost surely (a.s., for short)}.
$$

Indeed, suppose that $Y(s)>0$ a.s. for a certain $s\in(a,b)$. 
By continuity, $Y(t)>c>0$ a.s. in a neighborhood of $s$, 
$a<s-r<s<s+r<b, r>0$. Consider the process $\eta$ such that 
$\eta(t)=0$ a.s. on $[a,s-r]\cup[s+r,b]$ and  $\eta(t)>0$ 
a.s. on $(s-r,s+r)$, and $\eta(t)=1$ a.s. on 
$\left(s-\dfrac{r}{2},s+\dfrac{r}{2}\right)$.
Then, $\int_a^b Y(t)\cdot \eta(t)dt\geq rc>0$ a.s. Consequently, 
$E\left[\int_a^b Y(t) \cdot \eta(t)dt\right] >0$, 
which completes the proof of our claim.
Taking into account this result, and the fact that $\eta$ is an arbitrary process, 
we deduce the desired stochastic fractional Euler--Lagrange equation:
$\partial_{2} L+{_t^s}{D}{_b^ \alpha}\left[ 
\partial_3 L\right]
+{_a^s}{D}{_t^ \alpha}\left[ 
\partial_4 L\right]=0$.
The proof is complete.
\end{proof}

By adopting the same method as in the proof of Theorem~\ref{thm:SFE-Leq} 
and using our result of integration by parts for 
stochastic Caputo fractional derivatives, i.e.,
item (iii) of Lemma~\ref{lemmaIbyP}, we obtain 
the appropriate stochastic Caputo fractional Euler--Lagrange 
necessary optimality condition.

\begin{Theorem}[The stochastic Caputo fractional Euler--Lagrange equation]
\label{thm:SFE-Leq:C}
If $J \in C^{1}(H \times \mathbb{R} 
\times \mathbb{R} \rightarrow \mathbb{R})$ and $X\in C^1(I \rightarrow H)$
is an $F$-adapted stochastic process on $[a,b]$ with $E(X(t))\in AC([a,b])$ 
that is a minimizer~of 
\begin{equation*}
J[X]=E\left(\int_a^b 
L\left( t,X(t),{_a^{sC}}{D}{_t^ \alpha} X(t),
{_t^{sC}}{D}{_b^ \alpha} X(t)\right)dt\right)
\end{equation*}
subject to the fixed end points
$E(X(a))=X_a$ and $E(X(b))=X_b$,
then $X$ satisfies the following stochastic fractional Euler--Lagrange equation:
$$
\frac{\partial L}{\partial X}+{_t^{sC}}{D}{_b^ \alpha}\left[ 
\dfrac{\partial L}{\partial{_a^{sC}}{D}{_t^ \alpha}}\right]
+{_a^{sC}}{D}{_t^ \alpha}\left[ 
\dfrac{\partial L}{\partial{_t^{sC}}{D}{_b^ \alpha}}\right]=0.
$$
\end{Theorem}

\begin{Remark}
Note that the conclusions of Theorems~\ref{thm:SFE-Leq} 
and \ref{thm:SFE-Leq:C} are not contradictory: one conclusion
is valid for Riemann--Liouville derivative problems, 
while the other holds true for Caputo-type problems. 
The conclusions are proved in a similar manner by remarking 
that the additional quantity with parentheses, in the integration 
by parts theorem linked to the Caputo approach, 
vanishes under the condition that $X$ and $X^*$ verify 
the same initial and final conditions.
Note also that the assumptions of Theorems~\ref{thm:SFE-Leq} 
and \ref{thm:SFE-Leq:C} are necessary for the 
existence of left and right stochastic 
Riemann--Liouville/Caputo fractional derivative operators.
\end{Remark}


Our Theorems~\ref{thm:SFE-Leq} and \ref{thm:SFE-Leq:C}
give an extension of the Euler--Lagrange 
equations of the classical calculus of variations \cite{MR2004181}, 
stochastic calculus of variations \cite{R3}, 
and fractional calculus of variations \cite{R2}. 


\section{Examples}
\label{sec:5}

The best way to illustrate a new theory is by choosing simple examples.
We give two illustrative examples of the stochastic 
Riemann--Liouville fractional Euler--Lagrange equation
proved in Section~\ref{sec:4}: the first one inspired from quantum mechanics; 
the second chosen to allow a simple numerical solution to the obtained stochastic 
Riemann--Liouville fractional Euler--Lagrange equation.

\begin{Example}
\label{ex:01}
Let us consider the stochastic fractional variational problem
\eqref{eq:F}--\eqref{eq:BC} with
$$ 
L\left(t,X(t),{_a^s}{D}{_t ^ \alpha}X(t),  
{_t ^s}{D}{_b^ \alpha}X(t)\right)
=\frac{1}{2}\left(\frac{1}{2}m
\mid {_a^s}{D}{_t ^ \alpha}X(t) \mid ^2
+\frac{1}{2}m\mid {_t ^s}{D}{_b^ \alpha}X(t)\mid ^2\right)-V(X(t)), 
$$ 
where $X$ is a stochastic process of second order with 
$E(X(t)) \in AC([a,b])$ and $V$ maps $C'(I \rightarrow H)$ to $\mathbb{R}$. 
Note~that 
$$
\frac{1}{2}\left(\frac{1}{2}m\mid {_a^s}{D}{_t^ \alpha}X(t) \mid ^2
+\frac{1}{2}m \mid {_t^s}{D}{_b^ \alpha}X(t)\mid ^2\right)
$$ 
can be viewed as a generalized kinetic energy in 
the quantum mechanics framework. By applying our 
Theorem~\ref{thm:SFE-Leq} to the current variational 
problem, we get 
\begin{equation}
\label{eq:Newton:dyn:law}
\frac{1}{2}m \left[ {_a^s}{D}{_t^ \alpha}({_t^s}{D}{_b^ \alpha}X(t))
+{_t^s}{D}{_b^ \alpha}({_a^s}{D}{_t^ \alpha}X(t))\right]=grad V(X(t)),
\end{equation}
where $grad V$ is the gradient of $V$, which in this case
means the derivative of the potential energy of the system.
We observe that if $\alpha$ tends to zero and $X$ is a deterministic function, 
then relation \eqref{eq:Newton:dyn:law} becomes what is known in the physics 
literature as Newton's dynamical law: $m \ddot{X}(t)=grad V(X(t))$.
\end{Example}

The calculus of variations can assist us both analytically and numerically.
Now we give a numerical example, carried out with the help of 
the MATLAB computing environment \cite{MR3677002}.

\begin{Example}
\label{ex:02}
Let $\alpha := 0.25$, $a := 0.01$, $b := 0.99$,
$X_a := 1.00$, and $X_b := 1.00$. Consider the following 
variational problem \eqref{eq:F}--\eqref{eq:BC}:
\begin{gather*}
J[X]=\int_a^b {_a^s}{D}{_t^\alpha}X(t)\times {_t^s}{D}{_b^\alpha}X(t) \, dt
\longrightarrow \min,\\
E(X(a))=X_a, \quad E(X(b))=X_b,
\end{gather*}
where $X\in C^1(I\rightarrow H)$ with $E(X(t)) \in AC$ 
and ${_a^s}{D}{_\cdot ^ \alpha}X$ and ${_\cdot ^s}{D}{_b^ \alpha}X $ 
denote, respectively, the left and the right stochastic 
fractional Riemann--Liouville derivatives of order $\alpha$. 
Resorting again to Theorem~\ref{thm:SFE-Leq}, 
we obtain the following stochastic fractional Euler--Lagrange 
differential equation:
$$
{_a^s}{D}{_t^ {2\alpha}}X(t) + {_t^s}{D}{_b^ {2\alpha}}X(t)=0.
$$
Following \cite{R1}, we observe that 
${_a^s}{D}{_t^ \alpha}X(t)$ and $ {_t^s}{D}{_b^ \alpha}X(t)$ 
can be approximated as follows: 
$$
{_a^s}{D}{_t^ \alpha}X(t)={_a}{D}{_t^ \alpha}E(X(t))
\simeq \sum_{k=0}^N\frac{(-1)^{(k-1)}\alpha (E(X(t)))^{(k)}}{k!(k-\alpha)
\gamma(1-\alpha)}(t-a)^{(k-\alpha)}
$$ 
and
$$
{_t^s}{D}{_b^ \alpha}X(t)={_t}{D}{_b^ \alpha}E(X(t))
\simeq \sum_{k=0}^N\frac{-\alpha (E(X(t)))^{(k)}}{k!(k-\alpha)
\gamma(1-\alpha)}(b-t)^{(k-\alpha)}.
$$ 
Choosing $N=1$, we get the curve for $E(X(t))$ 
as shown in Figure~\ref{fig:01}.
\begin{figure}[H]
\centering
\includegraphics[scale=0.45]{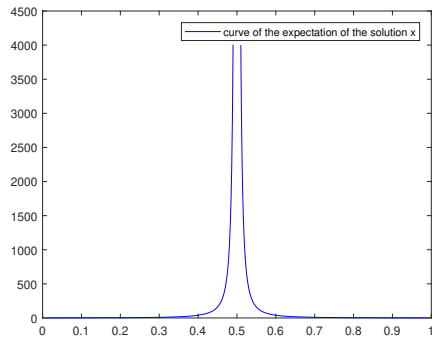}
\caption{Expectation of the extremal to the stochastic fractional 
problem of the calculus of variations of Example~\ref{ex:02}.}
\label{fig:01}
\end{figure}
One can increase the value of $N$ under the condition one adds 
a sufficient number of initial values related to some degrees 
of derivatives of $E(X(t))$. This particular question is similar 
to the standard fractional calculus and we refer 
the interested reader to the book \cite{R1}.
\end{Example}


\section{Conclusions}
\label{sec:6}

Numerous works related to the calculus of variations, addressing 
different optimization problems by means of classical, stochastic, 
and fractional derivatives through appropriate Euler-Lagrange equations, 
exist in the literature. To extend available results to 
a stochastic-fractional framework, we have established 
in this work new definitions associated to left and right stochastic 
Riemann--Liouville/Caputo fractional integrals and derivatives, together 
with some properties of boundedness, linearity, additivity and 
interaction between involved operators. Furthermore, we have proven 
new integration by parts theorems, according to the novel definitions, 
which have a central role for the establishment of the stochastic 
Riemann--Liouville/Caputo fractional Euler--Lagrange equations. 
The obtained stochastic Riemann--Liouville/Caputo fractional Euler--Lagrange 
equations generalize those available on the literature of fractional calculus.
Moreover, the results of the paper motivate readers and researchers to go on
and further develop  the theory now initiated. 

It is important to note that the mathematical background used 
in the original fractional calculus differs from what 
we have established here for the stochastic fractional case. 
Additionally, the six constructed definitions for the left and right 
stochastic Riemann--Liouville/Caputo integral/derivative operators, 
as well as proved integration by parts formulas, differ totally 
to those available in the fractional calculus theory: the first 
are applied to second order stochastic processes, and the second act 
on deterministic absolute continuous functions. Furthermore, our 
stochastic fractional Euler--Lagrange equations serve as necessary optimality
conditions to optimization problems subject to unknown stochastic processes  
that can be effectively approximated by numerical methods: such~equations 
might be transformed to ones subject to unknown  deterministic functions 
that are the expectation $E(X(t))$, for~instance, when the random variables  
$X(t)$ follow the assumption of normality, which~is instructive to 
approximate its expectation via stochastic fractional Euler--Lagrange equations 
determined statistically by the hypothesis test in inferential statistics.

We claim that the new  mathematical concepts we have introduced here
are more natural than those already available in the literature,
since it is intuitive and convenient to proceed 
via application of the expectation.

\vspace{6pt}
\authorcontributions{The authors equally contributed 
to this paper, read and approved the final manuscript. 
All authors have read and agreed to the published version of the manuscript.}

\funding{This research was funded by the Portuguese 
Foundation for Science and Technology (FCT),
grant~number UIDB/04106/2020 (CIDMA).} 
	

\acknowledgments{The authors are grateful to four anonymous reviewers 
for all their questions, comments and suggestions, which helped them 
to improve the clarity and quality of the paper.}

\conflictsofinterest{The authors declare no conflict of interest.} 


\reftitle{References}


\end{document}